\documentclass[12pt]{article}
\usepackage{hyperref}
\usepackage{amsfonts}
\usepackage{enumerate}
\begin{document}

\begin{center}
{\large\bf Homogeneous linear matrix difference equations of higher order: Singular case}

\vskip.20in
Charalambos P. Kontzalis$^{1}$ and\ Grigoris Kalogeropoulos$^{2}$\\[2mm]
{\footnotesize
$^{1}$Department of Informatics, Ionian University, Corfu, Greece\\
$^{2}$Department of Mathematics, University of Athens, Greece}
\end{center}

{\footnotesize
\noindent
\textbf{Abstract:} In this article, we study the singular case of an homogeneous generalized discrete time system with given initial conditions. We consider the matrix pencil singular and provide necessary and sufficient conditions for existence and uniqueness of solutions of the initial value problem.\\
\\[3pt]
\\[3pt]

\vskip.2in
\section{Introduction}
Many authors have studied generalized discrete \& continouus time systems, see [1-28], and their applications, see [29-35]. Many of these results have already been extended to systems of differential \& difference equations with fractional operators, see [36-45]. In this article, our purpose is to study the solutions of a generalized initial value problem of linear matrix difference equations into the mainstream of matrix pencil theory. Thus, we consider
\[
A_nX_{k+n}+A_{n-1}X_{k+n-1}+...+A_1X_{k+1}+A_0X_k=0_{m_1,1},
\]
with known initial conditions
\[
X_{k_0},X_{k_0+1},...,X_{k_0+n-1},
\]
where  $A_i, i=0,1,...,n \in \mathcal{M}({m_1 \times r_1;\mathcal{F}})$, (i.e. the
algebra of square matrices with elements in the field
$\mathcal{F}$) with  $X_k \in\mathcal{M}({m_1 \times
1;\mathcal{F}})$ and det$A_n$=0 if $A_n$ is square. In the sequel we adopt the following notations 
\[
    \begin{array}{c}
    Y_{k,1}=X_k,\\
    Y_{k,2}=X_{k+1},\\
    \dots \\
    Y_{k,n-1}=X_{k+n-2},\\
    Y_{k,n}=X_{k+n-1}.
    \end{array}
    \]
    and
   \[
   \begin{array}{c}
    Y_{k+1,1}=X_{k+1}=Y_{k,2},\\
    Y_{k+1,2}=X_{k+2}=Y_{k,3},\\
    \dots \\
    Y_{k+1,n-1}=X_{k+n-1}=Y_{k,n},\\
    A_nY_{k+1,n}=A_nX_{k+n}=-A_{n-1}Y_{k,n}-...-A_1Y_{k,2}-A_0Y_{k,1}.
    \end{array}
    \]
 Let $m_1n$=$r$ and $m_1n+r_1-m_1$=$m$. Then the above system can be written in Matrix form in the following way
\begin{equation}
FY_{k+1}=GY_k,
\end{equation}
with known initial conditions
\begin{equation}
Y_{k_0}.
\end{equation}      
Where  
\[
 F = \left[
\begin{array}{ccccc} 
I_{m_1}&0_{m_1, m_1}&...&0_{m_1, m_1}&0_{m_1, m_1}\\
0_{m_1, m_1}&I_{m_1}&...&0_{m_1, m_1}&0_{m_1, m_1}\\
\vdots&\vdots&\ddots&\vdots&\vdots\\
0_{m_1, m_1}&0_{m_1, m_1}&...&I_{m_1}&0_{m_1, m_1}\\
0_{m_1, r_1}&0_{m_1, r_1}&...&0_{m_1, r_1}&A_n
\end{array}
\right],
\]
\[
 G = \left[\begin{array}{cccc} 0_{m_1, m_1}&I_{m_1}&\ldots&0_{m_1, m_1}\\0_{m_1, m_1}&0_{m_1, m_1}&\ldots&0_{m_1, m_1}\\\vdots&\vdots&\ddots&\vdots\\0_{m_1, m_1}&0_{m_1, m_1}&\ldots&I_{m_1}\\-A_0&-A_1&\ldots&-A_{n-1}\end{array}\right].
\]
and 
\[
Y_k=\left[\begin{array}{cccc} Y_{k,1}^T Y_{k,2}^T \dots Y_{k,m_1}^T\end{array}\right]^T.
\]
With dimensions $F,G \in \mathcal{M}({r \times m;\mathcal{F}})$, (i.e. the
algebra of matrices with elements in the field
$\mathcal{F}$) with  $Y_k \in \mathcal{M}({m \times
1;\mathcal{F}})$. For the sake of simplicity we set
${\mathcal{M}}_m  = {\mathcal{M}}({m \times m;\mathcal{F}})$ and
${\mathcal{M}}_{rm}  = {\mathcal{M}}({r \times m;\mathcal{F}} )$. The matrices $F$ and $G$ can be non-square (when $r\neq m$) or square ($r = m$) and $F$ singular (det$F$=0).

\section{Singular matrix pencils: Mathematical background and notation}
In this section we will give the mathematical background and the notation that is used throughout the paper
\\\\
\textbf{Definition 2.1} Given $F,G\in \mathcal{M}_{rm}$  and an indeterminate $s\in\mathcal{F}$, the matrix
pencil $sF-G$ is called regular when  $r=n$ and  $\det(sF-G)\neq 0$.
In any other case, the pencil will be called singular.
\\\\
\textbf{Definition 2.2}
The pencil $sF-G$  is said to be \emph{strictly equivalent} to the
pencil $s\tilde F - \tilde G$ if and only if there exist nonsingular
$P\in\mathcal{M}_n$ and $Q\in\mathcal{M}_m$ such as
\[
    P({sF - G})Q = s\tilde F - \tilde G.
\]
In this article, we consider the case that the pencil is \emph{singular}. Unlike the case of the regular pencils, the characterization of a singular matrix pencil, apart from the set of the determinantal divisors requires the definition of additional sets of invariants, the minimal indices.  Let $\mathcal{N}_r$, $\mathcal{N}_l$ be right, left null space of a matrix respectively. Then the equations
\[
(sF-G)U(s)=0_{m, 1}
\]
and
\[
V^T(s)(sF-G)=0_{1, m}
\]
have solutions in $U(s), V(s)$, which are vectors in the rational vector spaces $\mathcal{N}_r(sF-G)$ and $\mathcal{N}_l(sF-G)$ respectively. The binary vectors X(s) and $Y^T(s)$ express dependence relationships among the colums or rows of $sF-G$ respectively. $U(s), V(s)$ are polynomial vectors. Let $d$=dim$\mathcal{N}_r(sF-G)$ and $t$=$\mathcal{N}_l(sF-G)$. It is known [46-53] that $\mathcal{N}_r(sF-G)$, $\mathcal{N}_l(sF-G)$, as rational vector spaces, are spanned by minimal polynomial bases of minimal degrees 
\[
\epsilon_1=\epsilon_2=...=\epsilon_g=0<\epsilon_{g+1}\leq...\leq\epsilon_d
\]
and 
\[
\zeta_1=\zeta_2=...=\zeta_h=0<\zeta_{h+1}\leq...\leq\zeta_t
\]
respectively. The set of minimal indices $\epsilon_i$ and $\zeta_j$ are known [46-53] as \textit{column minimal indices} (c.m.i.) and \textit{row minimal indices} (r.m.i) of sF-G respectively. 
To sum up in the case of a singular matrix pencil, we have invariants, a set of \emph{elementary divisors} (e.d.) and \emph{minimal indices}, of the following type:

\begin{itemize}
    \item e.d. of the type  $(s-a)^{p_j}$, \emph{finite elementary
    divisors} (nz. f.e.d.)
    \item e.d. of the type  $\hat{s}^q=\frac{1}{s^q}$, \emph{infinite elementary divisors}
    (i.e.d.).
    \item m.c.i. of the type $\epsilon_1=\epsilon_2=...=\epsilon_g=0<\epsilon_{g+1}\leq...\leq\epsilon_d$, \emph{minimal column indices}
    \item m.r.i. of the type $\zeta_1=\zeta_2=...=\zeta_h=0<\zeta_{h+1}\leq...\leq\zeta_t$, \emph{minimal row indices}
\end{itemize}
\textbf{Definition 2.3.} Let $B_1 ,B_2 ,\dots, B_n $ be elements of $\mathcal{M}_n$. The direct sum
of them denoted by $B_1  \oplus B_2  \oplus \dots \oplus B_n$ is
the blockdiag$\left[\begin{array}{cccc} B_1& B_2& \dots& B_n\end{array}\right]$.
\\\\
The existence of a complete set of invariants for singular pencils implies the existence of canonical form, known as Kronecker canonical form [46-53] defined by 
\[
sF_K  - Q_K :=
sI_p  - J_p  \oplus sH_q  - I_q \oplus sF_{\epsilon}-G_{\epsilon}\oplus sF_{\zeta}-G_{\zeta}\oplus 0_{h, g}
\]
where $sI_p  - J_p$ is uniquely defined by the set of f.e.d.
\[
  ({s - a_1 })^{p_1 } , \dots ,({s - a_\nu  }
 )^{p_\nu },\quad \sum_{j = 1}^\nu  {p_j  = p}
\]
of $sF-G$  and has the form
\[
    sI_p  - J_p  := sI_{p_1 }  - J_{p_1 } (
    {a_1 }) \oplus  \dots  \oplus sI_{p_\nu  }  - J_{p_\nu  }
    ({a_\nu  }) 
\]
The $q$  blocks of the second uniquely defined block $sH_q -I_q$ correspond to the i.e.d.
\[
  \hat s^{q_1} , \dots ,\hat s^{q_\sigma}, \quad \sum_{j =
  1}^\sigma  {q_j  = q}
\]
of $sF-G$  and has the form
\[
    sH_q  - I_q  := sH_{q_1 }  - I_{q_1 }  \oplus
    \dots  \oplus sH_{q_\sigma  }  - I_{q_\sigma}
\]
Thus, $H_q$  is a nilpotent element of $\mathcal{M}_n$  with index
$\tilde q = \max \{ {q_j :j = 1,2, \ldots ,\sigma } \}$, where
\[
    H^{\tilde q}_q=0_{q, q},
\]
and $I_{p_j } ,J_{p_j } ({a_j }),H_{q_j }$ are defined as
\[
   I_{p_j }  = \left[\begin{array}{ccccc} 
   1&0& \ldots & 0&0\\
   0& 1 &  \ldots&0 &0 \\
   \vdots & \vdots & \ddots & \vdots &\vdots \\
   0 & 0 & \ldots  & 0 &1
   \end{array}\right]
   \in {\mathcal{M}}_{p_j } , 
   \]
   \[
   J_{p_j } ({a_j }) =  \left[\begin{array}{ccccc}
   a_j  & 1 & \dots&0  & 0  \\
   0 & a_j  &   \dots&0  & 0  \\
    \vdots  &  \vdots  &  \ddots  &  \vdots  &  \vdots   \\
   0 & 0 &  \ldots& a_j& 1\\
   0 & 0 & \ldots& 0& a_j
   \end{array}\right] \in {\mathcal{M}}_{p_j }
\]
\[
 H_{q_j }  = \left[
\begin{array}{ccccc} 0&1&\ldots&0&0\\0&0&\ldots&0&0\\\vdots&\vdots&\ddots&\vdots&\vdots\\0&0&\ldots&0&1\\0&0&\ldots&0&0
\end{array}
\right] \in {\mathcal{M}}_{q_j }.
  \]
For algorithms about the computations of the Jordan matrices see [46-53]. For the rest of the diagonal blocks of $F_K$ and $G_K$, s$F_{\epsilon}-G_{\epsilon}$ and s$F_{\zeta}-G_{\zeta}$,  the matrices $F_{\epsilon}$, $G_{\epsilon}$ are defined as
\[ 
F_\epsilon=blockdiag\left\{L_{\epsilon_{g+1}}, L_{\epsilon_{g+2}}, ..., L_{\epsilon_d}\right\}
\]
 Where $L_\epsilon= \left[
\begin{array}{ccc} I_\epsilon & \vdots & 0_{\epsilon, 1}
\end{array}
\right]$, for $\epsilon=\epsilon_{g+1}, ..., \epsilon_d$
\[
G_\epsilon=blockdiag\left\{\bar L_{\epsilon_{g+1}}, \bar L_{ \epsilon_{g+2}}, ..., \bar L_{\epsilon_d}\right\}
\]  
Where $\bar L_\epsilon=\left[
\begin{array}{ccc} 0_{\epsilon, 1} & \vdots & I_\epsilon
\end{array}
\right]$, for $\epsilon=\epsilon_{g+1}, ..., \epsilon_d$. The matrices $F_{\zeta}$, $G_{\zeta}$ are defined as
\[ 
F_\zeta=blockdiag\left\{L_{\zeta_{h+1}}, L_{\zeta_{h+2}}, ..., L_{\zeta_t}\right\}
\]
 Where $L_\zeta= \left[
\begin{array}{c} I_\zeta \\ 0_{1, \zeta}
\end{array}
\right]$, for $\zeta=\zeta_{h+1}, ..., \zeta_t$
\[
G_\zeta=blockdiag\left\{\bar L_{\zeta_{h+1}}, \bar L_{ \zeta_{h+2}}, ..., \bar L_{\zeta_t}\right\}
\]
 Where $\bar L_\zeta=\left[
\begin{array}{c} 0_{1, \zeta}\\I_\zeta
\end{array}
\right]$, for $\zeta=\zeta_{h+1}, ..., \zeta_t$.

\section{Main results}

Following the given analysis in section 2, there exist non-singular matrices $P, Q$ such that 
\begin{equation}
PFQ=F_K, \quad PGQ=G_K.
\end{equation}
Let 
\begin{equation}
Q=\left[\begin{array}{ccccc}Q_p & Q_q &Q_\epsilon & Q_\zeta & Q_g\end{array}\right]
\end{equation}
where $Q_p\in \mathcal{M}_{rp}$, $Q_q\in \mathcal{M}_{rq}$, $Q_\epsilon\in \mathcal{M}_{r\epsilon}$, $Q_\zeta\in \mathcal{M}_{r\zeta}$ and $Q_g\in \mathcal{M}_{rg}$
\\\\
\textbf{Lemma 3.1.} System (1) is divided into five subsystems:
\begin{equation}
    Z_{k+1}^p = J_p Z_k^p
\end{equation}
the subsystem
\begin{equation}
    H_q Z^q_{k+1} = Z_k^q
\end{equation}
the subsystem
\begin{equation}
    F_\epsilon Z^\epsilon_{k+1}=G_\epsilon Z^\epsilon_k
\end{equation}
the subsystem
\begin{equation}
    F_\zeta Z^\zeta_{k+1}=G_\zeta Z^\zeta_k
\end{equation}
and the subsystem
\begin{equation}
    0_{h, g}\cdot Z^g_{k+1}=0_{h, g}\cdot Z^g_k
\end{equation}
\textbf{Proof.} Consider the transformation
\[
    Y_k=QZ_k
\]
By substituting this transformation into (1) we obtain
\[
    FQZ_{k+1}=GQZ_k.
\]
Whereby, multiplying by $P$ and using (3), we arrive at
\[
    F_KZ_{k+1}=G_K Z_k+PV_k.
\]
Moreover, we can write $Z_k$ as
\[
Z_k=\left[
\begin{array}{c} Z_k^p\\Z_k^q\\Z_k^\epsilon\\Z_k^\zeta\\Z_k^g
\end{array}
\right]
\]
where $Z_p^k\in \mathcal{M}_{p1}$, $Z_q^k \in \mathcal{M}_{q1}$, $Z_\epsilon^k\in \mathcal{M}_{\epsilon1}$, $Z_\zeta^k \in \mathcal{M}_{\zeta1}$ and $Z_g^k\in \mathcal{M}_{h1}$. Taking into account the above expressions, we arrive easily at the subsystems (5), (6), (7), (8), and (9).
\\\\
Solving the system (1) is equivalent to solving subsystems (5), (6), (7), (8), and (9).
\\\\
\textbf{Remark 3.1.} The subsystem (5) is a regular type system and its solution is given from, see [1-28].
\begin{equation}
Z^p_k=J_p^{k-k_0}Z^p_{k_0}.
\end{equation}
\textbf{Remark 3.2.} The subsystem (6) is a singular type system but its solution is very easy to compute, see [8-18].
\begin{equation}
Z^q_k=0_{q,1}.
\end{equation}
\textbf{Proposition 3.1.} The subsystem (7) has infinite solutions and can be taken arbitrary
\begin{equation}
Z_k^\epsilon=C_{k,1}.
\end{equation}
\textbf{Proof.} 
If we set 
\[
Z_k^\epsilon=\left[\begin{array}{c} Z_k^{\epsilon_{g+1}}\\Z_k^{\epsilon_{g+2}}\\\vdots\\Z_k^{\epsilon_d}\end{array}\right],
\]
by using the analysis in section 2, system (22) can be written as:
\[
blockdiag\left\{L_{\epsilon_{g+1}}, ..., L_{\epsilon_d}\right\}\left[\begin{array}{c} Z_{k+1}^{\epsilon_{g+1}}\\Z_{k+1}^{\epsilon_{g+2}}\\\vdots\\Z_{k+1}^{\epsilon_d}\end{array}\right]=blockdiag\left\{\bar L_{\epsilon_{g+1}}, ..., \bar L_{\epsilon_d}\right\}\left[\begin{array}{c} Z_k^{\epsilon_{g+1}}\\Z_k^{\epsilon_{g+2}}\\\vdots\\Z_k^{\epsilon_d}\end{array}\right].
\]
Then for the non-zero blocks a typical equation can be written as 
\[
\begin{array}{ccc} L_{\epsilon_i} Z_{k+1}^{\epsilon_i}=\bar L_{\epsilon_i} Z_k^{\epsilon_i} & , &i =g+1, g+2, ..., d, \end{array}
\]
or
\[
\left[\begin{array}{ccc} I_{\epsilon_i} & \vdots & 0_{{\epsilon_i}, 1}\end{array}\right]Z_{k+1}^{\epsilon_i}=\left[
\begin{array}{ccc} 0_{{\epsilon_i}, 1} & \vdots & I_{\epsilon_i}
\end{array}
\right]Z_k^{\epsilon_i},
\]
or
\[
\left[\begin{array}{ccccc} 1 &0&\ldots&0&0\\0&1&\ldots&0&0\\ \vdots &\vdots &\ldots&\vdots&\vdots\\0&0&\ldots&1&0\end{array}\right]\left[\begin{array}{c} z_{k+1}^{{\epsilon_i},1}\\ z_{k+1}^{{\epsilon_i},2}\\ \vdots \\ z_{k+1}^{{\epsilon_i},{\epsilon_i}} \\ z_{k+1}^{{\epsilon_i},{\epsilon_i}+1}\end{array}\right]=\left[\begin{array}{ccccc} 0 &1&\ldots&0&0\\0&0&\ldots&0&0\\ \vdots &\vdots &\ldots&\vdots&\vdots\\0&0&\ldots&0&1\end{array}\right]\left[\begin{array}{c} z_k^{{\epsilon_i},1}\\ z_k^{{\epsilon_i},2}\\ \vdots \\ z_k^{{\epsilon_i},{\epsilon_i}} \\ z_k^{{\epsilon_i},{\epsilon_i}+1}\end{array}\right],
\]
or
\[
\begin{array}{ccccc}  z_{k+1}^{{\epsilon_i},1}=z_k^{{\epsilon_i},2}\\ z_{k+1}^{{\epsilon_i},2}=z_k^{{\epsilon_i},3}\\\vdots\\ z_{k+1}^{{\epsilon_i},{\epsilon_i}}=z_{k+1}^{{\epsilon_i},{\epsilon_i}+1}.
\end{array}
\]
This is of a regular type system of difference equations with ${\epsilon_i}$ equations and ${\epsilon_i}+1$ unknowns. It is clear from the above analysis that in every one of the $d-g$ subsystems one of the coordinates of the solution has to be arbitrary by assigned total. The solution of the system can be assigned arbitrary
\[
Z_k^\epsilon=C_{k,1}
\]
\textbf{Proposition 3.2.} The subsystem (8) has the unique solution 
\begin{equation}
Z_k^\zeta=0_{\zeta,1}.
\end{equation}
\\\\
\textbf{Proof.} If we set
\[
Z_k^\zeta=\left[\begin{array}{c} Z_k^{\zeta_{h+1}}\\Z_k^{\zeta_{h+2}}\\\vdots\\Z_k^{\zeta_t}\end{array}\right],
\]
then the subsystem (8) can be written as:
\[
blockdiag\left\{L_{\zeta_{h+1}}, ..., L_{\zeta_t}\right\}\left[\begin{array}{c} Z_{k+1}^{\zeta_{h+1}}\\Z_{k+1}^{\zeta_{h+2}}\\\vdots\\Z_{k+1}^{\zeta_t}\end{array}\right]=blockdiag\left\{\bar L_{\zeta_{h+1}}, ..., \bar L_{\zeta_t}\right\}\left[\begin{array}{c} Z_k^{\zeta_{h+1}}\\Z_k^{\zeta_{h+2}}\\\vdots\\Z_k^{\zeta_t}\end{array}\right].
\]
Then for the non-zero blocks, a typical equation can be written as
\[
\begin{array}{ccc} L_{\zeta_j} Z_{k+1}^{\zeta_j}=\bar L_{\zeta_j} Z_k^{\zeta_j} & , & j=h+1, h+2, ..., t, \end{array}
\]
or
\[
\left[\begin{array}{c} I_{\zeta_j} \\ \cdots \\  0_{1, {\zeta_j}}\end{array}\right]Z_{k+1}^{\zeta_j}=\left[
\begin{array}{c} 0_{1,{\zeta_j}} \\ \cdots \\ I_{\zeta_j}
\end{array}
\right]Z_k^{\zeta_j},
\]
or
\[
\left[\begin{array}{cccc} 1 &0&\ldots&0\\0&1&\ldots&0\\ \vdots &\vdots &\ldots&\vdots\\0&0&\ldots&1\\0&0&\ldots&0\end{array}\right]\left[\begin{array}{c} z_{k+1}^{{\zeta_j},1}\\ z_{k+1}^{{\zeta_j},2}\\ \vdots \\ z_{k+1}^{{\zeta_j},{\zeta_j}}\end{array}\right]=\left[\begin{array}{cccc} 0 &0&\ldots&0\\1&0&\ldots&0\\ \vdots &\vdots &\ldots&\vdots\\0&0&\ldots&0\\0&0&\ldots&1\end{array}\right]\left[\begin{array}{c} z_k^{{\zeta_j},1}\\ z_k^{{\zeta_j},2}\\ \vdots \\ z_k^{{\zeta_j},{\zeta_j}}\end{array}\right],
\]
or
\[
\begin{array}{c}  z_{k+1}^{{\zeta_j},1}=0\\ z_{k+1}^{{\zeta_j},2}=z_k^{{\zeta_j},1}\\\vdots\\ z_{k+1}^{{\zeta_j},{\zeta_j}}=z_k^{{\zeta_j},{\zeta_j}-1}\\0=z_k^{{\zeta_j},{\zeta_j}}\end{array}.
\]
We have a system of ${\zeta_j}$+1 difference equations and ${\zeta_j}$ unknowns. Starting from the last equation we get the solutions
\[
\begin{array}{c}  z_k^{{\zeta_j},{\zeta_j}}=0\\ z_k^{{\zeta_j},{\zeta_j}-1}=0\\z_k^{{\zeta_j},{\zeta_j}-2}=0\\\vdots\\ z_k^{{\zeta_j},1}=0\end{array}
\]
Hence, the system (8) has the following unique solution
\[
Z_k^\zeta=0_{\zeta, 1}
\]
\textbf{Remark 3.3.}The subsystem (9) has an infinite number of solutions that can be taken arbitrary
\begin{equation}
Z_k^g=C_{k,2}
\end{equation}
We can state the following Theorem
\\\\
\textbf{Theorem 3.1.} Consider the system (1), with known initial conditions (2) and a singular matrix pencil $sF-G$. Then its solution is unique if and only if the c.m.i. are zero
\begin{equation}
dim\mathcal{N}_r(sF-G)=0
\end{equation}
and
\begin{equation}
Y_{k_0}\in colspan Q_p
\end{equation}
The unique solution is then given from the formula
\begin{equation}
    Y_k=Q_pJ_p^{k-k_0}Z_{k_0}^p
\end{equation}
where $Z_{k_0}^p$ is the unique solution of the algebraic system $Y_{k_0}=Q_pZ_{k_0}^p$. In any other case the system has infinite solutions.
\\\\
\textbf{Proof.} First we consider that the system has non zero c.m.i and non zero r.m.i. By using transformation 
\[
Y_k=QZ_k
\]
then from (10), (11), (12), (13) and (14) the solutions of the subsystems (5), (6), (7), (8) and (9) respectively are
\[
Z_k=\left[\begin{array}{c}
     Z_k^p  \\
     Z_k^q\\Z_k^\epsilon
     \\Z_k^\delta
     \\Z_k^g
     
    \end{array}\right] =\left[\begin{array}{c}
     J_p^{k-k_0}Z_{k_0}^p  \\
     0_{q, 1}\\C_{k,1}
     \\0_{t-h, 1}
     \\C_{k,2}
     
    \end{array}\right].
		\]
 Then by using (4)
\[
     Y_k = QZ_k =
     \left[\begin{array}{ccccc}Q_p & Q_q &Q_\epsilon & Q_\zeta & Q_g \end{array}\right]
     \left[\begin{array}{c}
     J_p^{k-k_0}Z_{k_0}^p  \\
     0_{q, 1}\\C_{k,1}
     \\0_{t-h, 1}
     \\C_{k,2}\end{array}\right] 
    \]
    and
    \[
    Y_k =
    Q_pJ_p^{k-k_0}Z_{k_0}^p+Q_\epsilon C_{k,1}+Q_g C_{k,2}
\]
Since $C_{k,1}$ and $C_{k,2}$ can be taken arbitrary, it is clear that the general singular discrete time system for every suitable defined initial condition has an infinite number of solutions. It is clear that the existence of c.m.i. is the reason that the systems (7) and consequently (9) exist. These systems have always infinite solutions. Thus a necessary condition for the system to have unique solution is not to have any c.m.i. which is equal to 
\[
dim\mathcal{N}_r(sF-G)=0.
\]
In this case the Kronecker canonical form of the pencil $sF-G$ has the following form
\[
sF_K  - Q_K :=
sI_p  - J_p  \oplus sH_q  - I_q \oplus sF_{\zeta}-G_{\zeta}
\]  
and then the system (1) is divided into the three subsystems (5), (6), (8) with solutions (10), (11), (13) respectively. Thus
\[
     Y_k = QZ_k =
     \left[\begin{array}{ccc}Q_p & Q_q & Q_\zeta\end{array}\right]
     \left[\begin{array}{c}
     J_p^{k-k_0}Z_{k_0}^p  \\
     0_{q, 1}
     \\0_{t-h, 1}
     \end{array}\right] 
    \]
		and
    \[
Y_k =
    Q_pJ_p^{k-k_0}Z_{k_0}^p.
\]
The solution that exists if and only if
\[
Y_{k_0}=Q_pZ_{k_0}^p,
\]
 or
\[
Y_{k_0}\in colspan Q_p.
\]
In this case the system has the unique solution 
\[
   Y_k=Q_pJ_p^{k-k_0}Z_{k_0}^p.
\]

\end{document}